\documentclass{article}

\usepackage{amssymb}
\usepackage{amsmath}

\usepackage[dvips]{graphicx}

\begin{document}

 
\begin{center}
\textbf{A REMARK ON CONTINUITY OF POSITIVE LINEAR FUNCTIONALS ON SEPARABLE BANACH $\ast$-ALGEBRAS}                                                                                                                                                                                     
\end{center}

\noindent \textbf{}

\begin{center}
M. El Azhari
\end{center}

\noindent \textbf{ } 

\noindent \textbf{Abstract.} Using a variation of the Murphy-Varopoulos Theorem, we give a new proof of the following R. J. Loy Theorem: Let $A$ be a separable Banach $\ast$-algebra with center $Z$ such that $ZA$ has at most countable codimension, then every positive linear functional on $A$ is continuous.

\noindent \textbf{}


\noindent \textbf{}


\noindent \textbf{} 
 
\noindent \textbf{}If $A$ is a $\ast$-algebra, a linear functional $f$ on $A$ is called positive if $f(x^{\ast}x)\geq 0$ for all $x\in A.$ A linear functional $g$ on $A$ dominates a linear functional $h$ on $A$ if $g-h$ is positive. Given an algebra $A$ and vector subspaces $X,Y$ of $A,\ XY $ will denote the vector subspace of $A$ spanned by the products $xy$ for $x\in X,\ y\in Y.$ Let $n\geq 2,\ X^{n}$ denotes the vector subspace of  $A$ spanned by the products $x_{1}\cdots x_{n}$ for $x_{i}\in X\ (1\leq i\leq n).$ A Hausdorff topological space $S$ is called a Souslin space if there is a complete separable metric space $P$ and a continuous mapping of $P$ onto $S$.

\noindent \textbf{}

\noindent \textbf{}We need two theorems and a preliminary proposition.

\noindent \textbf{}

\noindent \textbf{Theorem 1}([1, Theorem 5.5]).  Let $A$ be a complete separable metrizable topological vector space, $B$ a vector subspace of $A.$ If $B$ is a Souslin space and has at most countable codimension, then $B$ is closed and of finite codimension.

\noindent \textbf{}

\noindent \textbf{}In [5, Theorem], Varopoulos proved that if $A$ is a commutative Banach $\ast$-algebra, with continuous involution, such that $A^{3}$ is closed and of finite codimension, then every positive linear functional on $A$ is continuous. In [3, Corollary and Remark], Murphy gave another proof of this result without the assumption of continuity of the involution. In [4, Theorem 13.7], by the same methods in [3], Sinclair obtained the following improvement: Let $A$ be a Banach $\ast$-algebra with center $Z.$ If $Z^{2}A$ is closed and of finite codimension, then every positive linear functional on $A$ is continuous. Here we show that the above result works for $ZA^{2}.$

\noindent \textbf{}

\noindent \textbf{Theorem 2.}(Variation of the Murphy-Varopoulos Theorem). Let $A$ be a Banach $\ast$-algebra with center $Z.$ If $ZA^{2}$ is closed and of finite codimension, then every positive linear functional on $A$ is continuous.

\noindent \textbf{}

\noindent \textbf{Proof}: By hypothesis and [4, Lemma 13.6], it is sufficient to prove that every positive linear functional on $A,$ nonzero on $A^{2},$ dominates a continuous positive linear functional on $A,$ nonzero on $A^{2}.$ Let $f$ be a positive linear functional on $A,$ nonzero on $A^{2}.$ By the Schwarz inequality, we have $
\vert f(z^{\ast}xy)\vert^{2}\leq f_{z}(xx^{\ast})f(y^{\ast}y)$ for all $z\in Z,\ x\in A$ and $y\in A,$ where  $f_{z}(a)=f(z^{\ast}az)$ for all $a\in A.$ if $f_{z}$ is zero on $A^{2}$ for all $z\in Z,$ then $f$ is zero on $ZA^{2}$ and hence $f$ is continuous on $A.$ If there is $z\in A$ such that $f_{z}$ is nonzero on $A^{2},$ we can suppose $\Vert z^{\ast}z\Vert < 1.$ By the square root lemma [4, Lemma 13.1], there is $u\in A$ such that $u^{\ast}=u$ and $2u-u^{2}=z^{\ast}z.$ Let $x\in A, (f-f_{z})(x^{\ast}x)=f(x^{\ast}x-z^{\ast}x^{\ast}xz)=f(x^{\ast}x-x^{\ast}z^{\ast}zx)=f((x-ux)^{\ast}(x-ux))\geq 0.$ Thus $f$ dominates $f_{z}$ which is a continuous positive linear functional on $A$ by [4, Corollary 13.3].  

\noindent \textbf{}

\noindent \textbf{Proposition 3.} Let A be an algebra with center $Z.$ The following assertions are equivalent:
\begin{enumerate}
\item  $ZA^{m}$ has at most countable codimension for some $m\geq 1.$
\item  $ZA^{n}$ has at most countable codimension for all $n\geq 1.$
\end{enumerate}

\noindent \textbf{}

\noindent \textbf{Proof}: (1)$\Rightarrow$(2): We prove the implication by induction on $n.\ A^{2}$ and $ZA$ have at most countable codimension since $ZA^{m}\subset ZA\subset A^{2}.$ Suppose that $ZA^{n}$ has at most countable codimension. We have $A=ZA^{n}+E,\ A=A^{2}+F,\ E$ and $F$ are vector subspaces of $A$ with at most countable dimension. $A=A^{2}+F=(ZA^{n}+E)^{2}+F\subset ZA^{n+1}+E^{2}+F,$ then $A=ZA^{n+1}+E^{2}+F.$ So $ZA^{n+1}$ has at most countable codimension because $E^{2}+F$ has at most countable dimension.

\noindent \textbf{}

\noindent \textbf{Proof of R. J. Loy Theorem}([2, Theorem 2.1]): By Proposition 3, $ZA^{2}$ has at most countable codimension. Consider the continuous mapping $G:Z\times A\times A\to A,\  G(z,x,y)=zxy.$ Since $ZA^{2}$ is the linear span of $G(Z\times A\times A),$ it follows that $ZA^{2}$ is a Souslin space. By Theorem 1, $ZA^{2}$ is closed and of finite codimension, and so we can apply Theorem 2.

\noindent \textbf{}

\noindent \textbf{Remark}. The interest of the above proof is the observation that the argument for continuity in the Loy Theorem, which uses a result [2, Theorem 1.2] concerning multilinear mappings on products of separable Banach algebras, can be replaced by a variation of the well known theorem of Murphy-Varopoulos.

\noindent \textbf{}

\noindent \textbf{References}

\noindent \textbf{}

\noindent \textbf{}[1] J. P. R. Christensen, Topology and Borel structure, North-Holland Mathematical Studies, 10; Notas de Matematica 1974.

\noindent \textbf{}[2] R. J. Loy, Multilinear mappings and Banach algebras, J. London Math. Soc., 14 (1976), 423-429.

\noindent \textbf{}[3] I. S. Murphy, Continuity of positive linear functionals on Banach $\ast$-algebras, Bull. London Math. Soc., 1 (1969), 171-173.
 
\noindent \textbf{}[4] A. M. Sinclair, Automatic continuity of linear operators, London Math. Soc., Lecture Note Series 21, 1976.   

\noindent \textbf{}[5] N. Th. Varopoulos, Continuit\'{e} des formes lin\'{e}aires positives sur une  alg\`{e}bre  de Banach avec involution, C. R. Acad. Sci. Paris, 258 (1964), 1121-1124.

\noindent \textbf{}

\noindent \textbf{} 

\noindent \textbf{} Ecole Normale Sup\'{e}rieure

\noindent \textbf{} Avenue Oued Akreuch

\noindent \textbf{} Takaddoum, BP 5118, Rabat

\noindent \textbf{} Morocco
 
\noindent \textbf{} 

\noindent \textbf{} E-mail:  mohammed.elazhari@yahoo.fr

\end{document}